\documentclass[12pt,a4paper,leqno]{article}
\usepackage{etoolbox} \usepackage{comment}
\usepackage{amsfonts}
\usepackage[margin=1in]{geometry}
\usepackage{times}
\usepackage{amsthm,amssymb}
\usepackage{amsmath,xypic}
\usepackage{pdflscape}
\usepackage[export]{adjustbox} 
\usepackage[cal=boondoxo,calscaled=.96,scr=rsfs]{mathalpha} \usepackage{multicol}
\usepackage{hyphenat}
\usepackage[usenames,dvipsnames]{color}
\usepackage{rotating}
\usepackage{longtable}
\usepackage{caption}
\usepackage{tcolorbox}

\usepackage{epigraph}
\usepackage{enumitem}
\setlist[enumerate]{listparindent=0.5in}
\newcommand{\be}{\begin{equation}}
\newcommand{\ee}{\end{equation}}
\newcommand{\bes}{\begin{equation*}}
\newcommand{\ees}{\end{equation*}}
\newcommand{\bea}{\begin{eqnarray}}
\newcommand{\eea}{\end{eqnarray}}
\newcommand{\beas}{\begin{eqnarray}}
\newcommand{\eeas}{\end{eqnarray}}
\newcommand{\ben}{\begin{note}}
\newcommand{\een}{\end{note}}
\newcommand{\bexl}{\vskip0.1em\noindent\hrulefill\vskip1em\begin{ExerciseList}}
\newcommand{\eexl}{\end{ExerciseList}\hrulefill}

\newcommand{\bthm}{\begin{theorem}}
\newcommand{\ethm}{\end{theorem}}
\newcommand{\bpro}{\begin{prop}}
\newcommand{\epro}{\end{prop}}
\newcommand{\bcor}{\begin{corollary}}
\newcommand{\ecor}{\end{corollary}}
\newcommand{\bcon}{\begin{conjecture}}
\newcommand{\econ}{\end{conjecture}}
\newcommand{\bp}{\begin{proof}}
\newcommand{\ep}{\end{proof}}
\newcommand{\blem}{\begin{lemma}}
\newcommand{\elem}{\end{lemma}}
\newcommand{\bn}{\begin{note}}
\newcommand{\en}{\end{note}}
\newcommand{\benum}{\begin{enumerate}}
\newcommand{\eenum}{\end{enumerate}}
\newcommand{\bed}{\begin{defn}}
\newcommand{\eed}{\end{defn}}
\newcommand{\brem}{\begin{remark}}
\newcommand{\erem}{\end{remark}}

\newcommand{\btik}{\begin{tikzpicture}\begin{axis}[scale=0.5,axis y line=center, axis x line=middle]}
\newcommand{\etik}{\end{axis}\end{tikzpicture}}

\let\into=\hookrightarrow
\let\mapsto=\longmapsto

\newcommand{\upperRomannumeral}[1]{\uppercase\expandafter{\romannumeral#1}}

 \usepackage{tikz}
\usetikzlibrary{mindmap,backgrounds}
\usetikzlibrary{automata}
\usepackage{graphicx}
\usepackage[pagewise]{lineno}
\usepackage{stackengine}
\usepackage{stmaryrd}

\newcommand{\authornamemode}{\usepackage[backend=biber,style=authoryear]{biblatex}\setlength\bibitemsep{\baselineskip}}
\usepackage{titlesec}
\newtoggle{arxiv}
\toggletrue{arxiv}
\iftoggle{arxiv}{
	\usepackage{url}
\usepackage{natbib}
	\bibliographystyle{plainnat}
	\let\cite=\citep
	\typeout{In the arxiv mode}
	\usepackage{fancyhdr}
	\usepackage[colorlinks,citecolor=blue]{hyperref}
}
{
    \usepackage[pagewise]{lineno}
\authornamemode
	\usepackage{fancyhdr}
	\usepackage[colorlinks,citecolor=blue,colorlinks=true,hyperindex, citecolor=blue, urlcolor=blue]{hyperref}
	\let\cite=\parencite
\addbibresource{../../master/master6.bib}
\addbibresource{hoshi-bib.bib}
\addbibresource{mochizuki-bib.bib}
\addbibresource{uchida-bib.bib}
\addbibresource{mochizuki-flowchart.bib}
}
\newtoggle{draft}
\togglefalse{draft}

\newtheorem{theorem}[equation]{Theorem}      \newtheorem{lemma}[equation]{Lemma}          \newtheorem{corollary}[equation]{Corollary}  \newtheorem{proposition}[equation]{Proposition}

\theoremstyle{definition}
\newtheorem{conj}[equation]{Conjecture}

\theoremstyle{definition}
\newtheorem{defn}[equation]{Definition}
\theoremstyle{remark}

\theoremstyle{definition}
\newtheorem{remark}[equation]{Remark}

\newcommand{\para}{\subsection{}}
\newcommand{\parat}[1]{\subsection{#1}}
\titleformat{\subsection}[runin]{\normalfont\bfseries}{\thesubsection}{.5em}{}[{\ \ }]
\titlespacing{\subsection}{0pt}{1.5ex plus .1ex minus .2ex}{0pt}

\let\into=\hookrightarrow
\let\isom=\simeq

\newcommand{\A}{\mathscr{A}}
\newcommand{\abs}[1]{\left\vert#1\right\vert}

\newcommand{\bF}{{\bar{F}}}
\newcommand{\bQ}{{\bar{\Q}}}

\newcommand{\C}{{\mathbb C}}

\newcommand{\F}{{\mathbb F}}

\newcommand{\N}{\mathscr{N}}

\newcommand{\Q}{{\mathbb Q}}

\newcommand{\Z}{{\mathbb Z}}

\renewcommand{\O}{{\mathscr O}}
\renewcommand{\P}{{\mathbb P}}
\renewcommand{\wp}{{\mathfrak p}}

\newcommand{\fm}{{\mathfrak{M}}}

\newcommand{\invlim}{\varprojlim}
\let\fm=\fa
\newcommand{\mapright}[1]{{\xymatrix{{}\ar[r]^{#1}&{}}}}

\renewcommand{\bpro}{\begin{proposition}}
	\renewcommand{\epro}{\end{proposition}}
\renewcommand{\bcon}{\begin{conj}}
	\renewcommand{\econ}{\end{conj}}

\newcommand{\preliminary}{{\\ \normalsize \textcolor{blue}{Preliminary version for comments}}{\relax}}
\title{Comments on Arithmetic Teichmuller Spaces
	 \preliminary
}
\author{Kirti Joshi}

\newcommand{\Address}{\bigskip\noindent{\footnotesize\textsc{{Math. department, University of Arizona, 617 N Santa Rita, Tucson
		85721-0089, USA.}}\par\nopagebreak 
\noindent\textit{Email:}	\texttt{kirti@math.arizona.edu}}}

\begin{document}
	\maketitle
	\setcounter{tocdepth}{1}

\newcommand{\iut}{\cite{mochizuki-iut1, mochizuki-iut2, mochizuki-iut3,mochizuki-iut4}}
\newcommand{\topics}{\cite{mochizuki-topics1,mochizuki-topics2,mochizuki-topics3}}

\begin{abstract}
I provide some comments on Arithmetic Teichmuller Theory constructed in \cite{joshi-teich} and clarify its relationship to \iut.
\end{abstract}

\lhead{}

\iftoggle{draft}{\pagewiselinenumbers}{\relax}
\newcommand{\act}{\curvearrowright}
\newcommand{\lmp}{{\Pi\act\Ot}}
\newcommand{\lmpi}{{\lmp}_{\int}}
\newcommand{\lmpf}{\lmp_F}
\newcommand{\Om}{\O^{\times\mu}}
\newcommand{\Omf}{\O^{\times\mu}_{\bF}}
\renewcommand{\N}{\mathbb{N}}
\newcommand{\yoga}{Yoga}
\newcommand{\gl}[1]{{\rm GL}(#1)}
\newcommand{\bK}{\overline{K}}
\newcommand{\reptrip}{\rho:G_K\to\gl V}
\newcommand{\reptripp}[1]{\rho\circ\alpha:G_{\ifstrempty{#1}{K}{{#1}}}\to\gl V}
\newcommand{\benumlab}{\begin{enumerate}[label={{\bf(\arabic{*})}}]}
\newcommand{\ord}{\mathop{\rm ord}\nolimits}	
\newcommand{\kcs}{K^\circledast}
\newcommand{\lcs}{L^\circledast}
\renewcommand{\A}{\mathbb{A}}
\newcommand{\bfq}{\bar{\mathbb{F}}_q}
\newcommand{\tripod}{\P^1-\{0,1728,\infty\}}

\newcommand{\vseq}[2]{{#1}_1,\ldots,{#1}_{#2}}
\newcommand{\anab}[4]{\left({#1},\{#3 \}\right)\anabelmap\left({#2},\{#4 \}\right)}

\newcommand{\gln}{{\rm GL}_n}
\newcommand{\glo}[1]{{\rm GL}_1(#1)}
\newcommand{\glt}[1]{{\rm GL_2}(#1)}

\newcommand{\linv}{\mathfrak{L}}
\newcommand{\bedef}{\begin{defn}}
\newcommand{\eedef}{\end{defn}}
\renewcommand{\act}[1][]{\overset{#1}{\curvearrowright}}
\newcommand{\bfx}{\overline{F(X)}}
\newcommand{\anabelmap}{\leftrightsquigarrow}
\newcommand{\ban}[1][G]{\mathscr{B}({#1})}
\newcommand{\pit}{\Pi^{temp}}
 
 \newcommand{\bL}{\overline{L}}
 \newcommand{\bkm}{\bK_M}
 \newcommand{\vbk}{v_{\bK}}
 \newcommand{\vbkm}{v_{\bkm}}
\newcommand{\ocs}{\O^\circledast}
\newcommand{\ot}{\O^\triangleright}
\newcommand{\ocsk}{\ocs_K}
\newcommand{\otk}{\ot_K}
\newcommand{\ok}{\O_K}
\newcommand{\oko}{\O_K^1}
\newcommand{\oks}{\ok^*}
\newcommand{\Qpb}{\overline{\Q}_p}
\newcommand{\Qpbh}{\widehat{\overline{\Q}}_p}
\newcommand{\tr}{\triangleright}
\newcommand{\ocpt}{\O_{\C_p}^\tr}
\newcommand{\ocpf}{\O_{\C_p}^\flat}
\newcommand{\sG}{\mathscr{G}}
\newcommand{\sxfe}{\mathscr{X}_{F,E}}
\newcommand{\sxfep}{\mathscr{X}_{F,E'}}
\newcommand{\loglt}{\log_{\sG}}
\newcommand{\fc}{\mathfrak{t}}
\newcommand{\ku}{K_u}
\newcommand{\kup}{\ku'}
\newcommand{\kt}{\tilde{K}}
\newcommand{\sGpf}{\sG(\O_K)^{pf}}
\newcommand{\hgm}{\widehat{\mathbb{G}}_m}
\newcommand{\bE}{\overline{E}}
\newcommand{\sY}{\mathscr{Y}}
\newcommand{\syfe}{\mathscr{Y}_{F,E}}
\newcommand{\syfqp}[1]{\mathscr{Y}_{\cptl{#1},\Q_p}}
\newcommand{\syfqpe}[1]{\mathscr{Y}_{{#1},E}}
\newcommand{\fJ}{\mathfrak{J}}
\newcommand{\fM}{\mathfrak{M}}
\newcommand{\locvar}{local arithmetic-geometric anabelian variation of fundamental group of $X/E$ at $\wp$}
\newcommand{\fjxep}{\fJ(X,E,\wp)}
\newcommand{\fjxe}{\fJ(X,E)}
\newcommand{\fpc}[1]{\widehat{{\overline{\F_p(({#1}))}}}}
\newcommand{\cpt}{\C_p^\flat}
\newcommand{\cptl}[1]{\C_{p,{#1}}^\flat}
\newcommand{\fja}[1]{\fJ^{\rm arith}({#1})}
\newcommand{\ainfe}{A_{\inf,E}(\O_F)}
\renewcommand{\ainfe}{W_{\O_E}(\O_F)}
\newcommand{\gmh}{\widehat{\mathbb{G}}_m}
\newcommand{\sE}{\mathscr{E}}
\newcommand{\bpi}{B^{\varphi=\pi}}
\newcommand{\bpip}{B^{\varphi=p}}
\newcommand{\onto}{\twoheadrightarrow}

\newcommand{\cpmax}{{\C_p^{\rm max}}}
\newcommand{\xan}{X^{an}}
\newcommand{\yan}{Y^{an}}
\newcommand{\bPi}{\overline{\Pi}}
\newcommand{\bPit}{\bPi^{\rm{\scriptscriptstyle temp}}}
\newcommand{\Pit}{\Pi^{\rm{\scriptscriptstyle temp}}}
\renewcommand{\pit}[1]{\Pi^{\scriptscriptstyle temp}_{#1}}
\newcommand{\pitk}[2]{\Pi^{\scriptscriptstyle temp}_{#1;#2}}
\newcommand{\pio}[1]{\pi_1({#1})}
\newcommand{\fTeich}{\widetilde{\fJ(X/L)}}
\newcommand{\ssep}{\S\,} \newcommand{\vphi}{\varphi}
\newcommand{\sgt}{\widetilde{\sG}}
\newcommand{\sxqp}{\mathscr{X}_{\cpt,\Q_p}}

\setcounter{tocdepth}{1}

\tableofcontents

\togglefalse{draft}
\newcommand{\FF}{\cite{fargues-fontaine}}
\iftoggle{draft}{\pagewiselinenumbers}{\relax}

\newcommand{\attportion}{Sections~\ref{se:number-field-case}, \ref{se:construct-att}, \ref{se:relation-to-iut}, \ref{se:self-similarity} and \ref{se:applications-elliptic}}

\newcommand{\Pib}{\overline{\Pi}}
\newcommand{\four}{Sections~\ref{se:grothendieck-conj}, \ref{se:untilts-of-Pi}, and \ref{se:riemann-surfaces}}
\section{Two anabelian approaches to $p$-adic and adelic Teichmuller Theory}
I  offer some comments on \cite{joshi-teich} and its relationship to \iut. Especially I discuss the similarities and differences between the two theories. This document is intended to be accessible to a broad mathematical audience. To keep the length of this document short, I have avoided the technical aspects of the two theories and  I will not recall full statements of the theorems regarding Arithmetic Teichmuller Spaces  established in \cite{joshi-teich} and \cite{joshi-teich-estimates} but restrict myself to providing an impressionistic view of the core ideas. In Section~\ref{appendix-intereminacies}, I provide a discussion of how Mochizuki's Indeterminacies appear from the point of view of \cite{joshi-teich} and in  \ref{se:classical-case-example}, I discuss an illustrative example which will help readers arrive at a better understanding of some aspects of the two theories.

 I would like to acknowledge one mathematician who provided comments and suggestions on several versions of this manuscript but who wishes to remain anonymous. Another mathematician (anonymous by request) recently suggested some stylistic improvements for greater readability and is also acknowledged here.

\para Before proceeding, it will be useful for  the readers to keep the following remark in mind:

The absolute Grothendieck conjecture for hyperbolic curves over number fields (known by \cite{mochizuki96-gconj} and \cite{tamagawa97-gconj}) does not in any way  preclude the existence of the classical Teichmuller spaces at any archimedean prime (of a number field). This is because the classical Teichmuller space deals with the underlying complex analytic space $\xan/\C$ and this analytic structure  can be deformed  unconstrained by the absolute Grothendieck conjecture over number fields (as evidenced by the existence of classical Teichmuller spaces).

\para As far as I understand, Mochizuki's main point in \iut\ is that a similar phenomenon occurs at any non-archimdean prime: there is a genuine $p$-adic arithmetic Teichmuller  space at every non-archimedean prime unconstrained by the absolute Grothendieck conjecture for hyperbolic curves over $p$-adic fields.  This point is, perhaps, not stated in \iut\ with the sufficient clarity and emphasis it deserves  and to be sure such a space is not constructed explicitly in \iut\ (to my knowledge) \emph{but} Mochizuki discovered and described most of its key properties. These are encoded in the basic toolkit of \iut--Hodge theaters, Frobenioids and their complicated calculus described in \iut. 

\para In \cite{joshi-teich}, pursing an independent approach, I have arrived at the existence of $p$-adic arithmetic Teichmuller spaces using the theory of algebraically closed perfectoid fields \cite{scholze12-perfectoid-ihes}. The properties (\cite[\ssep 1.4]{joshi-teich}) of the arithmetic Teichmuller spaces I construct are quite similar to the ones discovered by Mochizuki  and have a natural geometric explanation in my approach  using \cite{fargues-fontaine}.

\para In a complete parallel to the theory of classical Teichmuller spaces, my approach  demonstrates very clearly that the Berkovich analytic structure of $\xan/\C_p$ deforms unconstrained by the validity of the absolute Grothendieck conjecture for hyperbolic curves over $p$-adic fields (\cite[Theorem 3.9.1 and Theorem 3.15.1]{joshi-teich}). These deformations arise from the fact the algebraically closed perfectoid field $\C_p$ has non-trivial deformations as a valued field \cite{kedlaya18} (note that Berkovich analytic space datum is dependent on the valued base field used in the construction of the space) and their anabelian consequences appear in my theory (\cite[Theorem 3.9.1]{joshi-teich}) proved using the principle of permanence of the tempered fundamental groups (of geometrically connected, smooth quasi-projective varieties) under extension of algebraically closed, complete valued fields due to \cite{lepage-thesis}. 

\para Mochizuki's approach in \iut\ may be termed ``group theoretic approach to $p$-adic arithmetic Teichmuller Theory'' while my approach \cite{joshi-teich} may be termed ``Berkovich analytic function theoretic approach  to $p$-adic arithmetic Teichmuller Theory'' and is closer to the more familiar complex function theoretic approach to classical Teichmuller Theory.  However let me also be perfectly clear that while my approach to the problem of constructing $p$-adic arithmetic Teichmuller spaces is fundamentally different from that of \iut,  it certainly could not have existed without the backdrop of Mochizuki's work. 

\para At any rate, regardless of how one arrives at it, the conclusion is inescapable: for each prime number $p$ and for each non-archimedean prime $\wp|p$ of a number field $L$, and for each geometrically connected, smooth, quasi-projective variety $X/L$ (Mochizuki works with hyperbolic curves but my results are broader--work in all dimensions) there is a fundamental and intrinsic degree of freedom which is also anabelian in nature in that the triple (here $L_\wp$ is the completion of $L$ at $\wp$ and $\C_p=\widehat{\bar{L}}_\wp$ is its completed algebraic closure)
\be
	\label{eq:fundamental-triple}\pi_1^{temp}(\xan/L_\wp)\hookleftarrow \pi_1^{temp}(\xan/\C_p)
\ee
consisting of the tempered fundamental group $\pi_1^{temp}(\xan/L_\wp)$ of the associated Berkovich analytic space  $\xan/L_\wp$  and the canonical inclusion  of the subgroup $\pi_1^{temp}(\xan/\C_p)$ corresponding to  $\xan/\C_p$ (i.e. the inclusion of the geometric tempered fundamental group) remains fixed up to isomorphism (this also implies that the \'etale fundamental group and the geometric \'etale fundamental group remain fixed up to isomorphism by \cite[Proposition~4.4.1]{andre03}) but the pair of Berkovich analytic spaces $(\xan/L_\wp,\xan/\C_p)$  associated to it deform unconstrained by any version (global or $p$-adic) of the absolute Grothendieck conjecture. Notably, since isomorphs of \eqref{eq:fundamental-triple} can arise from distinct geometric data, one cannot treat all the isomorphs of \eqref{eq:fundamental-triple} on an equal footing without ignoring the geometric data giving rise to them and Teichmuller Theory, classical  and $p$-adic theory of \cite{joshi-teich} (resp. \iut), is about the geometric data (resp. group theoretic data in \iut) giving rise to the isomorphs.

\para By its very construction, the theory of \cite{joshi-teich} is anabelian in nature (for aforementioned reasons) and hence it is fully compatible with the anabelian theory of \iut. Notably, starting with distinct geometric data $(\xan/L_\wp,\xan/\C_p)$ provided by  \cite{joshi-teich}  and passing (if required) to the triple \eqref{eq:fundamental-triple} of group theoretic data arising from them and applying the  constructions of \iut\ of  (local) Frobenioids, Hodge-Theaters  (which have components which are assembled from local Frobenioids) to each of the geometric data, one obtains corresponding objects arising from geometrically distinct local Berkovich data (at each prime $p$ including $p=\infty$). Since topologically distinct geometric data give rise to these Frobenioids, Hodge-Theaters etc., so in the theory of \cite{joshi-teich}, these objects of \iut\ must be considered distinct--and hence, as far as I understand, also in \iut.

\para Another important assertion of \iut\ is that there exists a fundamental and non-trivial invariance of valuation scale  in his theory and this is an operation incompatible with ring structures and should be seen as relating two distinct scheme theories in a some non-ring theoretic way. This assertion of \iut\ that there is some $p$-adic Teichmuller Theory equipped with some sort of scale invariance was already quite interesting for me.\footnote{Let me remind the readers that theories with scale invariance built into them are quite interesting--classical theory of Electromagnetism (with no sources), and fractals or self-similar objects are  examples of scale  invariant theories and objects. }

\para In \cite{joshi-teich} I have independently demonstrated this valuation scaling property (for the theory of  \cite{joshi-teich}) and my (new) observation is that this invariance of scale arises from the manifestly non-ring/scheme theoretic operation (permissible and  performed) on the set of closed points of degree one of the complete  Fargues-Fontaine curve $\sxqp$ constructed in \cite{fargues-fontaine}, and given explicitly by $([a]-p)\mapsto ([a^r]-p)$ (for any $a\in\O_{\cpt}, 0< \abs{a}<1$,  and any $0<r$ in the value group of $\cpt$). Notably the principle ideals $([a]-p)$ and $([a^r]-p)$ of $W(\O_{\cpt})$ provide closed points of the Fargues-Fontaine curve $\sxqp$ whose residue fields $K, K_r$ (say) are not necessarily  isomorphic as valued fields (by \cite{kedlaya18}) and hence this valuation scaling operation, in fact, alters  Berkovich geometries over these residue fields i.e. the pairs of analytic spaces $(X/E,\xan/K)$ and $(X/E,\xan/K_r)$  are not isomorphic while both provide \eqref{eq:fundamental-triple} i.e. the fundamental group and the geometric fundamental subgroup provided by these pairs are isomorphic \cite{joshi-teich}. A more precise description of this valuation scaling property of \cite{joshi-teich} is: Berkovich geometry altering valuation scaling property. In particular, Berkovich geometry altering valuation scaling property--and more generally the existence of distinct underlying Berkovich function theories of \cite{joshi-teich}, \cite{joshi-teich-estimates} is a demonstrably more robust notion (as is proved in \cite{joshi-teich}) than the view of \iut, that its operations must be seen as ``gluing together models of conventional scheme theories.'' 

\para A fancier version of this valuation scaling property is Mochizuki's Theta-link \iut, which says that Mochizuki's theory admits not only valuation scaling property for a single valuation but also admits  valuation scaling property for many valuations simultaneously. I establish this simultaneous valuation scaling property explicitly  in \cite{joshi-teich-estimates} building on  \cite{joshi-teich} by constructing an explicit, uncountable subset $\Sigma_F$ (with the said simultaneous valuation scaling property) of the self-product $\sxqp^{\ell^*}$ (where $\ell^*=\frac{\ell-1}{2}$ for  an odd prime $\ell$) of the complete Fargues-Fontaine curve $\sxqp$. Since I arrived at this set while trying to understand Mochizuki's Theta-link critical to \iut, I call this subset \emph{Mochizuki's Ansatz}. By construction Mochizuki's Ansatz has a rich geometry and exists independently of any anabelian considerations. 

\para To understand the significance of this set $\Sigma_F$, fix a canonically punctured Tate elliptic curve $X/L_\wp$ over a finite extension of $\Q_p$. The existence of this set  permits one to ``collate'' (in the precise sense of \cite{joshi-teich}) values (at some chosen $\ell$-torsion points) of a chosen Theta  function on $X/L_\wp$ (considered as a function on $\xan/\C_p$)  as the data of Berkovich analytic spaces $(\xan/L_\wp,\xan/\C_p)$ moves  while \eqref{eq:fundamental-triple} stays fixed. This leads to the construction of the Theta-values set $\tilde{\Theta}\subset B_{\Q_p}^{\ell^*}$ (where  $B_{\Q_p}$ is the Fargues-Fontaine ring for the datum $F=\cpt,E=\Q_p$) of \cite{joshi-teich-estimates}. Using the rich structure of Mochizuki's Ansatz, one can also  prove a lower bound for the size of $\tilde{\Theta}$ analogous to the one proved in \cite{mochizuki-iut3}. My constructions   take place in the theory of \cite{joshi-teich} but introduce a number of innovations (such as the introduction of the set $\Sigma_F$, working in the value of group of $\cpt$ which eliminates the need to track identifications of value groups in \iut) which leads to a conceptually cleaner proof of  (\cite[Theorem 10.1.]{joshi-teich-estimates} which is the analog of \cite[Corollary 3.12]{mochizuki-iut3} in the theory of \cite{joshi-teich} etc.) and   are certainly independent of \iut. Nevertheless these constructions are very broadly  guided by Mochizuki's rubric.

\section{Differences between the two approaches} Note that no Diophantine inequalities are claimed in \cite{joshi-teich} but there are other differences between \cite{joshi-teich} and \iut\  and let  me point them out now.

\para An important point to note is this: in \iut, Mochizuki works globally i.e. over a number field from the very beginning, and this is one of the important causes of confusion--as it gives the impression that global rigidity forced by the validity of the absolute Grothendieck conjecture (over number fields) precludes any local variation. In fact one must first understand that at each local picture (i.e. at each prime of the number field)  one encounters fluidity of geometric data unconstrained by any version (local or global) of the absolute Grothendieck conjecture (\cite{joshi-teich}). Then one assembles the locally fluid geometric data into an adelic one by general adelic methods. The approach of \cite{joshi-teich} develops the local picture in detail first and then proceeds to the adelic picture.    So the approach of \cite{joshi-teich} is local to adelic. This is a  major difference between how the two approaches proceed.

\para The next  most obvious, and major difference between \cite{joshi-teich} and \iut\ is the use of the theory of algebraically closed perfectoid fields and Berkovich spaces over such fields in \cite{joshi-teich}. In \iut\ Mochizuki makes extensive use of Anabelian Reconstruction Theory which, as far as I understand, cannot reconstruct the algebraically closed, perfectoid field $\C_p$ (as its anabelian reconstruction is obstructed by a well-understood result of \cite{mochizuki-local-gro}). The point of \cite{joshi-teich} is that the field $\C_p$ itself deforms as a valued field and hence deforms Berkovich geometry along with it. Thus approach of  \cite{joshi-teich} avoids Anabelian Reconstruction Theory altogether by working with geometric data from the very beginning and so geometric data does not have to be reconstructed at all. 

\emph{Personally}, I have found that Anabelian Reconstruction Theory is quite difficult to work with and my understanding is that it is also the one of important source of technical difficulties in  \iut. 

Let me remark that in \cite{joshi-anabelomorphy}, I have advocated the view that the notion of \emph{anabelomorphism of schemes} (an anabelomorphism of schemes is an isomorphism between the fundamental groups of the schemes) is an equivalence relation on schemes and also the notion of \emph{amphoric quantities} associated to schemes (i.e quantities associated to the scheme which are invariants of its anabelomorphism class)  as  mathematically more robust and general notions. This allows us to view anabelomorphic schemes and  quantities associated with them through an optik which is natural, geometric and also broader  than the lens of anabelian reconstruction theory of \topics\ (anabelian reconstruction theory is available only in dimension one at the moment). Hopefully readers will be persuaded (in the long run) that the view point of \cite{joshi-teich} provides a conceptually cleaner construction of arithmetic Teichmuller spaces which is both anabelian and geometric in nature.

\renewcommand{\fm}{\mathfrak{m}}

\para Many assertions of \iut\ are couched in the language of \emph{indeterminacies}. This becomes necessary if one does not admit the existence of $p$-adic arithmetic Teichmuller spaces. Following classical analogy will be useful: if one is presented with an isomorph of the topological fundamental group of a Riemann surface, then the  datum of the  point of the classical Teichmuller space giving rise to the isomorph  is necessarily indeterminate. I came to  understand indeterminacies of \iut\ in this sense. 

In contrast, the language of indeterminacies is not needed in \cite{joshi-teich} because one works with   geometric data i.e. pairs of Berkovich analytic spaces giving rise to \eqref{eq:fundamental-triple} from the very beginning. If one wants to use the language of indeterminacies in the context of \cite{joshi-teich} then one can say that the pair of Berkovich spaces $(\xan/L_\wp,\xan/\C_p)$ giving rise to \eqref{eq:fundamental-triple} is an indeterminate and thus admits this particular indeterminacy. This fundamental ``indeterminacy'' of \cite{joshi-teich} i.e. the existence of many pairs of analytic spaces $(\xan/L_\wp,\xan/\C_p)$ giving rise to the triple \eqref{eq:fundamental-triple} corresponds  to the three fundamental indeterminacies  Ind1, Ind2, Ind3 of \cite[Page 12]{mochizuki-iut3} (at one prime $\wp$ in this discussion). This is explicated for the readers in \ref{appendix-intereminacies} of this paper--\emph{a reading of it  requires a bit more technical material from \cite{joshi-teich} and so it may be skipped for an initial reading of this paper}. So the fact that fixing \eqref{eq:fundamental-triple} does not fix the pair of analytic spaces $(\xan/L_\wp,\xan/\C_p)$  in the theory of \cite{joshi-teich} leads  to indeterminacies which are related to the anabelian indeterminacies of  \iut.

\para In \iut, Mochizuki's construction of his ``theta-values set'' required for the assertion of \cite[Corollary 3.12]{mochizuki-iut3} takes place in $\log$-shell (or self-products of log-shells). This is simply the image of the natural continuous homomorphism $\log(\O_{L_\wp}^*) \to L_\wp$ (Mochizuki has shown that this is an amphoric structure in the sense of \cite{joshi-anabelomorphy}). Working with $\log$-shells arising from different isomorphs of local Frobenioids  requires a careful book-keeping of log-shells carried out in \iut. 

\para In contrast, my computations (\cite{joshi-teich-estimates}) are carried out in the Fargues-Fontaine ring $B$ (more precisely in self-products of this ring). There are many advantages for working with the ring $B$ as opposed to Mochizuki's log-shells, notably all the theta values arising from different pairs of analytic spaces corresponding to \eqref{eq:fundamental-triple} can be seen in one location simultaneously. Second $B$ is a Fr\'echet algebra (so one can make estimates with respect to the Fr\'echet structure) and it is equipped with a Frobenius $\vphi:B\to B$ and a continuous action of $G_{\Q_p}$ (so one can pass to ring of invariants), and thirdly one has a natural (Fr\'echet) subspace $ B^{\vphi=p}\subset B$ (where the action of $\vphi$ is multiplication by $p$) and by \cite[Proposition 4.1.3 and Example 4.4.7]{fargues-fontaine} one has a continuous isomorphism $1+\fm_{\O_{\cpt}}\mapright{x\mapsto \log([x])} B^{\vphi=p}\subset B$. So $B^{\vphi=p}$ is a log-shell in the sense of Mochizuki contained in $B$! From the perspective of \cite{joshi-teich} and \cite{joshi-teich-estimates} it is also possible to work in $B^{\vphi=p}$ (instead of $B$) but that is not carried out in the present version of \cite{joshi-teich-estimates}.

\para Let me record an important point: $p$-adic arithmetic Teichmuller Theory of \cite{joshi-teich} is indexed by the pairs $(X/L_\wp,\xan/\C_p)$ such that \eqref{eq:fundamental-triple} remains fixed up to isomorphism, while as far as I understand the indexing set of \iut\ is arbitrary. \emph{Nevertheless   I show that my theory (\cite{joshi-teich} and \cite{joshi-teich-estimates}) has all the key properties of the theory developed in \iut.}

\para Finally let me come to the main reason why $p$-adic arithmetic Teichmuller spaces (regardless of whether one prefers group theoretic version or the more geometric version) are of fundamental interest. The proof of the geometric Szpiro inequality due to \cite{szpiro79}, \cite{kim97},  \cite{bogomolov00}, \cite{zhang01},  \cite{beauville02} and others are based on the existence of archimedean Teichmuller spaces. Mochizuki's idea in \iut\ is that his discovery of genuine $p$-adic arithmetic Teichmuller spaces can be utilized to provide a proof of the Szpiro inequality in the number setting as well. Especially the proof in \iut\ is closest to the group theoretic proofs of the geometric Szpiro inequality given by \cite{bogomolov00} and \cite{zhang01}. See \ref{se:classical-case-example} for discussion of a classical Teichmuller Theory example which readers my find useful in this context.

\para For a more precise discussion of the relationship between the theory of \cite{joshi-teich} and \iut, readers are referred to \cite{joshi-teich}, \cite{joshi-teich-estimates} and readers may also find my papers \cite{joshi-anabelomorphy}, \cite{joshi-formal-groups} and \cite{joshi-gconj} useful.

\section{Mochizuki's Three Indeterminacies from the perspective of \cite{joshi-teich}}\label{appendix-intereminacies}
In this section I want to discuss how the fact that \eqref{eq:fundamental-triple} does not determine pairs of analytic spaces $(\yan/E',\yan/K)$ providing \eqref{eq:fundamental-triple} gives rise to Mochizuki's Indeterminacies Ind1, Ind2, Ind3 (\cite[Page 12]{mochizuki-iut3}) which play a key role in \cite[Theorem 3.11]{mochizuki-iut3} and \cite[Corollary 3.12]{mochizuki-iut3}. Let me explicate this for the readers (\emph{this portion is a bit more technical and independent of \iut\ and so may be skipped for initial reading}).

\para Let me say that once one has constructed Arithmetic Teichmuller Spaces, the notion of Indeterminacies in the sense of Mochizuki's \iut, becomes largely redundant in the theory of \cite{joshi-teich}. \emph{Especially one should view Mochizuki's Indeterminacies as place-holders for a Teichmuller Space in \iut.} So it is useful to understand how the existence of Arithmetic Teichmuller Spaces of \cite{joshi-teich} provides Indeterminacies  essentially in Mochizuki's sense.

\para From the point of view of \cite{joshi-teich}, there is one and only one fundamental indeterminacy in the Arithmetic Teichmuller Theory of \cite{joshi-teich} namely the exact sequence \eqref{eq:fundamental-triple} does not uniquely determine the pairs of spaces $(Y/E',\yan/K)$ giving rise to \eqref{eq:fundamental-triple} and this fundamental indeterminacy gives rise to indeterminacies identifiable or corresponding with Mochizuki's Indeterminacies Ind1, Ind2, Ind3. Let me exlain how each arises in my theory.

\parat{Ind1}\label{Ind1} A pair $(\yan/E',\yan/K)$ gives rise to \eqref{eq:fundamental-triple} means by definition (\cite{joshi-teich}) that one has an anabelomorphism $\pit{X/L_{\wp}}\isom \pit{Y/E'}$ and hence by \cite[\ssep 8.7]{joshi-teich} (this is a consequence of an old result of Mochizuki) one has an isomorphism $G_{L_\wp}\isom G_{E'}$ i.e. $L_\wp$ and $E'$ are anabelomorphic $p$-adic fields (this terminology was introduced in \cite{joshi-anabelomorphy}). As is well-known (by an old results of \cite{yamagata76}, \cite{jarden79}), anabelomorphy of $p$-adic fields does not imply necessarily that $E'$ is  isomorphic to $L_\wp$. Thus if one fixes \eqref{eq:fundamental-triple}--the $p$-adic field $L_\wp$ is not determinable from \eqref{eq:fundamental-triple} and hence this is an indeterminacy. This is precisely Mochizuki's Indeterminacy Ind1 (at one prime $\wp$, in \iut\ Mochizuki works with all primes of $L$ simultaneously). 
	
\parat{Ind2} Let $\hgm/\Z_p$ be the multiplicative  formal group. Let me assume (for simplicity--and this is more than adequate for \iut\ but my results are more general) that the  perfectoid field $K$  (in the pair $(\yan/E',\yan/K)$ giving rise to \eqref{eq:fundamental-triple})  has its tilt $K^\flat\isom\cpt$ and hence one has a pair $(K\supset \Q_p,K^\flat\isom \cpt)$--such pairs are parameterized by the Fargues-Fontaine curve $\syfqp{}$ (by \cite{fargues-fontaine}). Let me mention a rather delicate point of \cite{joshi-teich} here: from the pair $(\yan/E',\yan/K)$ one also has  the pair $(K\supset E',K^\flat\isom\cpt)$, but this is a point of $\sY_{\cpt,E'}$, on the other hand  because of \ref{Ind1}  the field $E'$ is only determined up to an anabelomorphism with $E$. So  the non-trivial point here is that the Teichmuller space of \cite{joshi-teich} ``glues'' together all $\sY_{\cpt,E'}$ for $E'$ anabelomorphic to $E$. Working with $(K\supset\Q_p,K^\flat\isom \cpt)$, i.e. with $\syfqp{}$ (this forgets $E'$), one has the property that one has an isomorphism of topological groups (\cite[Theorem 8.29.1]{joshi-teich}) $$\invlim_{x\mapsto  x^p}\hgm(\O_K)= \invlim_{x\mapsto x^p}(1+\fm_{\O_K}) \isom \invlim_{x\mapsto x^p}(1+\fm_{\C_p}) =\invlim_{x\mapsto  x^p} \hgm(\O_{\C_p})$$
	that is, the isomorphism class of the topological group $$\widetilde{\hgm}(\O_{\C_p}):=\invlim_{x\mapsto  x^p}(1+\fm_{\C_p})$$ is determined  by the pairs $(Y/E',\yan/K), (X/E,\xan/\C_p)$ but \eqref{eq:fundamental-triple} does not determine these pairs (uniquely). 
	
	In \iut\ and especially in \cite{mochizuki-iut3}, Mochizuki works with the topological group
	$$\O_{\bQ_p}^{*\mu}=\O_{\bQ_p}^{*}/\mu({\Qpb})=(1+\fm_{\Qpb})/(\mu_{p^\infty}(\Qpb))$$ (here $\mu({\Qpb})$ resp. $\mu_{p^\infty}(\Qpb)$ is the group of roots of unity resp. group of $p$-power roots of unity in $\Qpb$) and arbitrary isomorphisms between two isomorphs of $\O_{\bQ_p}^{*\mu}$. \emph{Especially, and more precisely, in Mochizuki's theory--see \cite[Page 10]{mochizuki-iut2}, a local Frobenioid datum arising from $X/L_\wp$ at a prime $\wp|p$ determines the topological isomorphism class of the topological group $\O_{\bQ_p}^{*\mu}$.} 
	
	Notably   the topological group $\widetilde{\hgm}(\O_{\C_p})$ and arbitrary topological isomorphisms between two isomorphs of this topological group (such as the one exhibited by the above displayed equation) plays in \cite{joshi-teich} the role  played by the topological group $\O_{\bQ_p}^{*\mu}$ and arbitrary isomorphisms between two isomorphs of this group in \iut.
	
	So in \cite{joshi-teich}, the analog of Mochizuki's Indeterminacy Ind2 can again be more explicitly described arising from the indeterminacy of the pairs giving rise to \eqref{eq:fundamental-triple}. 
	\parat{Ind3} This indeterminacy arises in \cite{joshi-teich} as follows: let $(\yan/E',\yan/K)$ be a pair giving rise to \eqref{eq:fundamental-triple} and satisfying  $K^\flat=\cpt$. Now  one learns from \cite{matignon84}, \cite{kedlaya18} that there exist many isometric embeddings $\cpt\into \cpt$ which do not extend to isometries $\cpt\mapright{\isom}\cpt$. This says that $\cpt$ is a non-trivially self-similar,  algebraically closed, complete valued field of characteristic $p$. The analog of Mochizuki's indeterminacy Ind3 in \cite{joshi-teich} is this fundamental self-similarity property of $\cpt$. It is important to recognize this property of $\cpt$ as providing an indeterminacy: a pair of analytic spaces $(Y/E', \yan/K)$ with $K^\flat=\cpt$ giving rise to \eqref{eq:fundamental-triple} determines $K^\flat=\cpt$ but one has the freedom to pass to another copy of $\cpt$ by a strict isometric embedding! Consequences of this are detailed in \cite[\ssep 11]{joshi-teich}.

\section{An illustrative example}\label{se:classical-case-example} Let me provide an example which in my opinion provides a good  illustration of the core strategies of  \iut\ and \cite{joshi-teich}, \cite{joshi-teich-estimates}  to arrive at Diophantine applications of the existence of the respective approaches to Arithmetic Teichmuller Theory.  The core strategy is due to Mochizuki (\cite[Corollary 3.12]{mochizuki-iut3}) but the two implementations (\cite{mochizuki-iut3} and \cite{joshi-teich-estimates}) differ significantly.

\para Let us say we have some connected, hyperbolic Riemann surface $\Sigma=(\Sigma,g)$ equipped with its hyperbolic, conformal metric $g$ and a continuous path $\gamma:[0,1]\to \Sigma$. Suppose that one is interested in the length $\ell_{\Sigma}(\gamma)$ of this path. This length  of course  depends on the metric $g$ (hence the subscript $\Sigma$ in $\ell_{\Sigma}$). Now one is  interested in some non-trivial upper bound on $\ell_{\Sigma}(\gamma)$. One could do the following: take a subset $\Sigma\in\Theta\subset T_{\Sigma}$ of Riemann surfaces in the Teichmuller space $T_\Sigma$. Using the quasi-conformal mappings provided in the datum for $\Sigma'\in T_{\Sigma}$, one can view $\gamma$ as a continuous path in $\Sigma'$ (denoted again by $\gamma$). Then one has a trivial bound
$$L(\gamma,\Theta):={\rm Supremum}\left\{\ell_{\Sigma'}(\gamma):\Sigma'\in \Theta\right\}\geq \ell_{\Sigma}(\gamma).$$
If we can find a suitable subset $\Theta$ where the supremum $L(\gamma,\Theta)$ is easy to evaluate or to bound  from above then our problem of bounding $\ell_{\Sigma}(\gamma)$  is solved! \emph{The proof of the geometric Szpiro conjecture \cite{bogomolov00} and \cite{zhang01} in fact do (broadly) proceed along these lines!}

\para Let me make some observations about this example:
\benum[label={{\bf(\arabic{*})}}]
\item The key point is that there exists some Teichmuller space $T_\Sigma$ which provides many different Riemann surfaces homeomorphic to $\Sigma$.
\item One is free to chose any $\Theta\subset T_{\Sigma}$ so long as it contains $\Sigma$ and so long as  $L(\gamma,\Theta)$ can be bounded from above in some way. 
\item \emph{Importantly}: there is nothing which gets glued to anything! 
\item The only possible gluing (if one wants to call it by this name) is in the existence of quasi-conformal mappings $f:\Sigma'\to \Sigma$ (necessarily homeomorphisms) which allow us to move $\gamma$ from $\Sigma$ to a $\Sigma'\in \Theta$. 
\item \emph{Especially:} one can consider $\gamma$ as a path on a fixed topological surface $\Sigma$ equipped with many different complex structures provided by $\Theta$ (i.e. the surface and the path datum is fixed but the complex geometry of the surface is being deformed).
\item A better term to use here in place of `gluing' is collating of of the numerical data $\ell_{\Sigma'}(\gamma)$ (using $\Theta$)!
\item \emph{To obtain a non-trivial bound  we need  a set $\Theta$ containing surfaces $\Sigma'$ for which at least some of the lengths are strictly larger than $\ell_{\Sigma}(\gamma)$.}
\eenum

\para  One may view $\Sigma$ as a fixed topological surface with distinct quasi-conformal structures providing distinct elements of $\Theta\subset T_{\Sigma}$ (my view) or one may view the quasi-conformal structures on $\Sigma$ as providing  distinctly  labeled copies of the topological fundamental group of the surface $\Sigma$ (Mochizuki's view)!  
\newcommand{\uH}{\mathfrak{H}}
\newcommand{\sM}{\mathscr{M}}
\para Assume $\Sigma$ is compact and of genus one (i.e. $\Sigma$ is  of conformal type $(1,0)$). Then one can identify the Teichmuller space $T_{\Sigma}$ with the upper half-plane $T_{\Sigma}\isom\uH\subset \C$  and the coarse moduli $\sM_1$ of compact Riemann surfaces of genus one can be identified as its quotient (for the M\"obius action of ${\rm PSL}_2(\Z)$ on $\uH$) $$T_\Sigma\isom \uH\to \sM_1\isom \uH/{\rm PSL}_2(\Z)$$ and one can choose the subset $\Theta\ni\Sigma$ to be the fibre (in $T_{\Sigma}$) over $[\Sigma]\in \sM_1$.
One may also similarly describe the situation for other conformal types $(g,n)$. \emph{Notably} one can find $\Theta\ni\Sigma$ even if the isomorphism class of the Riemann surface $\Sigma$  is fixed.
\newcommand{\fS}{\mathfrak{S}}
\newcommand{\fSt}{\mathfrak{S}_{Teich}}

\para Now in \iut\ (and in \cite{joshi-teich-estimates}), one wants to do  this i.e. find such a suitable set $\Theta$ in the $p$-adic context (for every prime $p$) with $\ell_{\Sigma}(\gamma)$ replaced by the $p$-adic absolute value of the Tate parameter at each prime of bad semi-stable reduction. Mochizuki's construction of such a set $\Theta$ is described in \cite{mochizuki-iut3}. The construction of such a set $\Theta$ using the theory of \cite{joshi-teich} is carried out in \cite{joshi-teich-estimates}. 

\para More generally, in Diophantine geometry,  one is interested in heights of sections of line bundles evaluated at some algebraic points. The heights are dependent on local (Arakelov) metric data and the insight of \iut\ and especially of \cite{joshi-teich} is that the local $p$-adic metric data moves for  non-trivial geometric reasons and thus creates the possibility of arriving at estimates by considering some conveniently chosen (sub)set in the totality of metric data available (for each prime $p$) as a bounding set for the quantities of interest.

\para In the  $p$-adic context, and especially in \cite[Corollary 3.12]{mochizuki-iut3}, Mochizuki views his labels as providing distinct ``arithmetic holomorphic structures'' or even more precisely ``arithmetic holomorphic structures (with prescribed indeterminacies)'' (`arithmetic' as in arising from the nature of the coefficients) \cite[Page 20]{mochizuki-iut1}, \cite[Page 35]{mochizuki-iut2}.

In \cite{joshi-teich},  \emph{I substantially deepen this idea}: in my theory labels correspond to distinct (arithmetic) Berkovich analytic (i.e. holomorphic) structures.
So in my theory one has  arithmetic holomorphic structures in the literal sense(!) and in fact I demonstrate   (in \cite{joshi-teich}) that there are many distinct (arithmetic) Berkovich analytic structures. Hence there are also many distinct labels in Mochizuki's sense too.

\para In Mochizuki's case, labels (and their calculus) form the Teichmuller Theory in the title of \iut\ and in my case (\cite{joshi-teich}) one has a Teichmuller Theory  in the classical sense describing variations of (Berkovich) holomorphic structures and equipped with a  similar calculus of labels.

\iftoggle{arxiv}{
\bibliography{hoshi-bib,mochizuki-bib,uchida-bib,mochizuki-flowchart,../../master/master6.bib}
}
{
	\printbibliography
}

\Address
\end{document}